# A STUDY OF GENERALIZED COVARIANT HAMILTON SYSTEMS ON GENERALIZED POISSON MANIFOLD

GEN WANG

Abstract. Since the basic theoretical framework of generalized Hamilton system is not perfect and complete, there are often some practical problems that can not be expressed by generalized Hamilton system. The generalized gradient operator is defined by the structure function on manifold to improve the basic theoretical framework of the whole generalized Hamilton system. The generalized structural Poisson bracket

$$\{f,g\} = \{f,g\}_{GPB} + G(s,f,g) = \{f,g\}_{GPB} + f\{s,g\}_{GPB} - g\{s,f\}_{GPB}$$

is defined as well on manifolds. The geometric bracket is also given, and the covariant extension form of the generalized Hamilton system directly related to the structure function, the generalized covariance, is further obtained–generalized covariant Hamilton system, It includes thorough generalized Hamiltonian system and S-dynamic system.

Contents

## 1. Introduction

Hamilton structure originates from classical mechanics, which makes some stability problems of classical mechanics better handled, and the theory of infinite-dimensional Hamilton system expands the scope of classical mechanics. Current research results show that Poisson structure in mechanical systems mainly comes from two aspects: first, from classical symmetry. Hamilton system is generated by reduction, such as Euler equation of fixed-point motion of free rigid body. Secondly, there exists Poisson structure in phase space of momentum mapping, such as fixed-point motion of heavy rigid body. For a given system, the establishment of Hamilton structure consists of two steps, i.e. writing Hamilton function and Poisson brackets [1]. The former can be composed of system. Conserved quantities are obtained. It can also be proved theoretically that Hamilton structure exists in any conservative system, but it is not easy to find practical Poisson bracket to make Hamilton structure simple. The phase space of mechanical systems often has a cotangent bundle structure, on which the symplectic form naturally defines a Poisson bracket. From this viewpoint, it is often an effective way to obtain a simpler Hamilton structure by using reduction theory [2, 3].

In 1953, Pauli [4], a Swiss physicist, discovered a generalization of non-canonical variables in Hamilton mechanics when he studied quantization of non-local field theory. In 1959, Martin [5], a British mathematician and physicist, tried to generalize the Hamilton method in order to apply it to systems without Lagrange functions. Their result, as a generalization theory of Hamilton system, is called generalized Hamilton system. Later, scientists used generalized Poisson bracket to define generalized Hamilton system, and more importantly. For the sake of simplicity and convenience, the further development of Hamilton mechanics has been promoted.

Generalized Hamilton system on $\mathbb{R}^r$ is a natural generalization of classical Hamilton system defined on $\mathbb{R}^{2n}$. The phase space of classical Hamilton system can only be even dimension, while the phase space of generalized Hamilton system can be arbitrary finite dimension or even infinite dimension. The phase space of generalized Hamilton system is a Poisson manifold, which usually has the characteristics of symplectic layer structure. Each symplectic layer is an invariant manifold of generalized Hamiltcm system. The generalized Hamilton system above is the classical Hamilton system on symplectic manifolds [6, 7].

Since 1950s, the theory of generalized Hamilton system has developed rapidly, but some practical problems can not be expressed by generalized Hamilton system. For example [8], the differential equation of fixed-point rotation of rigid body under external moment must add an additional term to the generalized Hamilton equation of Euler case, which becomes the generalized Hamilton system with additional terms is established.

Generalized Hamilton system dynamics plays an important role in many fields of science and technology. It has been widely valued in mathematics, mechanics, physics, engineering science and many other fields. It has made a series of important progress in theory, calculation and application. However, its basic theoretical framework is not perfect, and a lot of work still needs to be done. Further research is needed [6, 8, 9, 10].

The first question is: can a nonlinear system always be expressed as a generalized Hamilton system and how? Or under what conditions can it be expressed as a generalized Hamilton system? This is the so-called Hamiltonian realization. Generalized Hamilton implementation is a challenging research topic. Although there are some international research results on this issue, most of them are descriptive, which seriously hinders the application of generalized Hamilton implementation in nonlinear systems.

In this paper, the structural function determined by manifolds is introduced to discuss the more general theoretical framework of generalized Hamilton systems. The generalized structural Poisson bracket is used to define the generalized covariant Hamilton systems. The structural properties of generalized Hamilton systems on general manifolds are studied, and the dynamics of generalized Hamilton systems is improved. It can solve some practical problems. A perfect version of the generalized covariant Hamilton system is obtained. Thus, the main innovations of this paper are organized as follows:

(1) By using the covariant derivative operator to define a new general Poisson bracket called the generalized structural Poisson bracket (GSPB).
(2) Using new Poisson bracket to define a new complete covariant dynamics called generalized covariant Hamilton system (GCHS) that generalizes the GHS on manifold.

Generalized covariant Hamilton system answers the first question above. A non-linear system can not always be expressed as a generalized Hamilton system, but can always be expressed as a generalized covariant Hamilton system. Or a non-linear system can be expressed as a generalized covariant Hamilton system with consideration of structural functions. The generalized covariant Hamilton system is a complete theoretical system and systematic, which makes the application of the generalized covariant Hamilton system in non-linear systems eliminate some obstacles to the realization of the generalized Hamilton system.

Here is a brief description of the framework of the generalized structural Poisson bracket (**GSPB**) theory

**GPB** → **GSPB** → **GCHS** ⟵ **GHS** ← **GPB**

↙ ↘ ↓

**S-dynamics** **TGHS**

## 2. Canonical Poisson bracket and canonical Hamilton equation

The motion of an $n$ freedoms mechanical system can be described in a $2n$-dimensional phase space $P$ consisting of generalized coordinates and generalized conjugate momentum. For a conservative system, as long as the total energy $H : P \to \mathbb{R}$ of the system is known, the motion of the system can be described by $H$ as the following canonical Hamilton equation shown [11, 12].

$$\dot{q}_i = \frac{\partial H}{\partial p_i},\ \dot{p}_i = -\frac{\partial H}{\partial q_i},\ i = 1, 2, \cdots, n \tag{2.1}$$

The canonical equation can also be expressed in Poisson brackets on $P$:

$$\dot{q}_i = \{q_i, H\}_{CPB},\ \dot{p}_i = \{p_i, H\}_{CPB} \tag{2.2}$$

where $\{,\}_{CPB}$ is a canonical Poisson bracket on function space $C^r(p)$.

$$\{F_1, F_2\}_{CPB} = \sum_i \left( \frac{\partial F_1}{\partial q_i}\frac{\partial F_2}{\partial p_i} - \frac{\partial F_1}{\partial p_i}\frac{\partial F_2}{\partial q_i} \right),\ \forall F_j \in C^r(P) \tag{2.3}$$

Hamilton mechanics developed under the above framework is the most important component of classical mechanics. However, the Hamllton system defined above can only be defined in even-dimensional phase space. For odd-dimensional systems, it can not be described by the above methods, which restricts the application scope of Hamllton's viewpoint. In order to apply Hamilton's viewpoint to a wider range of dynamical systems, many famous scholars (S. Lie, a mathematician dating back to the last century) have made various generalizations, making it possible to define appropriate generalized Poisson structures on any finite dimensional manifold, which have all the main properties of (2.3). Later in the 1960s, V.I. Arnold studied the flow of incompressible ideal fluid from Hamllton's point of view, thus introducing Hamilton's point of view into the problem of infinite dimensional phase space and defining the generalized Poisson bracket structure in phase space. Then, according to Arnold's idea, J. Marsden et al. solved a series of stability problems of infinite dimensional dynamical systems by using the generalized Poisson bracket theory. The basic idea of the theory, which has been sorted out by the system, has become an important energy-Casimir method and energy-momentum method in the study of non-linear stability problems [12].

## 3. Generalized Poisson bracket (GPB) and generalized Hamilton system

For given functons $f, g \in C^\infty(M)$, generalized Poisson bracket (GPB) [1] is expressed as [1, 2, 3, 4, 5, 12]

$$\{f, g\}_{GPB} = \nabla^T f J \nabla g = \sum_{i,j} J_{ij}(x)\, \partial_i f \partial_j g \equiv J_{ij}(x)\, \partial_i f \partial_j g,\ i, j = 1, 2, \cdots, n \tag{3.1}$$

[1] PB: Poisson bracket or CPB: Classical Poisson bracket; GPB: Generalized Poisson bracket; GSPB: Generalized Structural Poisson Bracket; HS: Hamilton System; GHS: Generalized Hamilton System; TGHS: Thorough generalized Hamiltonian system; GCHS: Generalized Covariant Hamilton System.

It has the following important properties :

(1) Antisymmetry$\{f,g\}_{GPB} = -\{g,f\}_{GPB}$
(2) Bilinearity:$\{\lambda f + \mu g, h\}_{GPB} = \lambda\{f,h\}_{GPB} + \mu\{g,h\}_{GPB}$
(3) Jacobi identity: $\{f,\{g,h\}_{GPB}\}_{GPB} + \{g,\{h,f\}_{GPB}\}_{GPB} + \{h,\{f,g\}_{GPB}\}_{GPB} = 0$
(4) Leibnitz identity:$\{f\cdot g,h\}_{GPB} = f\cdot\{g,h\}_{GPB} + g\cdot\{f,h\}_{GPB}$
(5) non-degeneration: if for $f$,$\{f,g\}_{GPB} = 0$, then $g$ must be a constant.

where $f,g,h$ are elements in $C^\infty(P)$, $\lambda,\mu$ are arbitrary real number.

Generalized Poisson bracket (GPB) can be defined by the first four properties above.

**Definition 3.1.** [12] The generalized Poisson bracket on the smooth manifold $M$ is an operation on the smooth function space $C^\infty(M)$, for $f,g \in C^\infty(M)$, it exists $\{f,g\}_{GPB} \in C^\infty(M)$, the operation satisfies the condition $1-4$, $(M,\{,\}_{GPB})$ is correspondingly called Poisson manifolds.

The dimension of $M$ is not specified in this definition. It can be arbitrary finite dimension (especially odd dimension) or even infinite dimension.

Since GPB satisfies bilinear and Leibnitz identities, it actually defines derivative operations on $C^\infty(M)$. So for each smooth function $H \in C^\infty(M)$, it exists a vector field $X_H$ such that

$$X_H \cdot f = \{f,H\}_{GPB}$$

holds for all $f \in C^\infty(M)$.

**Definition 3.2.** [12] The vector field $X_H$ is called generalized Hamilon vector field determined by function $H$ on Poisson manifolds $M$, $H$ is Hamilton function.

Hence, generalized Poisson bracket $\{f,H\}_{GPB}$ is the derivative of $f$ in the direction of vector field $X_H$. Especially, if $(x_1,\cdots,x_m)$ is a local coordinates near point $p$ in $M$, then $\{x_i,H\}_{GPB}$ defines a derivative of coordinates $x_i$ along the direction of $X_H$, the integral curve of $X_H$ can then be parameterized with time $t$ such that

$$\dot{x}_i = \frac{dx_i}{dt} = \{x_i,H\}_{GPB}\,,\ i = 1,\cdots,m \tag{3.2}$$

This is the generalized Hamilton equation of motion. Therefore, for arbitrary function $f(x_1,\cdots,x_m) \in C^\infty(M)$, the rate of change along $X_H$ is

$$\dot{f} = \frac{df}{dt} = \{f,H\}_{GPB}$$

**Definition 3.3.** [6, 12] The skew-symmetric matrix $J = (J_{ij}(x))$ of order $m\times m$ is called structural matrix of the GPB on $M$, where $J_{ij}(x) = \{x_i,x_j\}_{GPB}$.

## 4. Generalized Hamilton System (GHS)

The expression of the generalized Hamilton system is as follows [1, 2, 6, 12]

$$\frac{dx}{dt} = J(x)\nabla H(x)\,,\ x \in \mathbb{R}^m$$

where $J(x)$ is structural matrix, $\nabla H(x)$ is gradient of Hamilton function.

**Proposition 4.1.** [12] *The necessary and sufficient condition for defining a GPB for the arbitrary function matrix $J = (J_{ij}(x))$ on $R^m$ is*

(1) *antisymmetry, $J_{ij}(x) = -J_{ji}(x)$*
(2) *Jacobi identity $J_{il}\frac{\partial J_{jk}}{\partial x_l} + J_{jl}\frac{\partial J_{ki}}{\partial x_l} + J_{kl}\frac{\partial J_{ij}}{\partial x_l} = 0$*

*where $i,j,k = 1,\ldots,n$.*

Note that condition 2 is actually Jacobi's identity. As Generalized Poisson bracket is defined as (3.1), by using GPB, the generalized Hamilton equation (3.2) can be further written as follows:

$$\dot{x}_i = \{x_i, H\}_{GPB} = J_{ij}\frac{\partial H}{\partial x_j},\ i = 1, \ldots, n \tag{4.1}$$

(1) when structural matrix $J(x) = J_0 = \begin{pmatrix} 0 & I_n \\ -I_n & 0 \end{pmatrix}$ holds, then the system is the classical Hamilton system (CHS).
(2) When the structural matrix element is a homogeneous linear function of $x$, namely, $J_{ij} = c_{ij}^k x_k$, the corresponding Poisson bracket is Lie-Poisson bracket.
(3) The dimension of classical Hamilton's structure matrix must be even dimension, and the generalized Hamilton can be arbitrary dimension.

The generalized Hamiltonian system is a mathematical description of the evolution of physical systems, which involves some basic concepts and theorems about time and phase space. Overall, the study of generalized Hamiltonian systems is a product of the intersection of physics, dynamics, and mathematical theory. The study of generalized Hamiltonian systems originated from the discovery of the "Three Laws of Mechanics". Since Newton, people have gradually discovered some laws under physical mechanisms. However, with the continuous development of science and technology, Newton's laws seem to be increasingly unable to adapt to the scientific research needs of this era, leading to the development of generalized mechanical systems based on Lagrange mechanics and Hamilton mechanics. In the promotion of generalized Hamiltonian system research, with the continuous development and cross integration of mathematics, physics, and computer science, more research directions have emerged. For example, generalized Hamiltonian systems in classical mechanics, quantum mechanics, and relativistic mechanics, with research scope including analogies of quantum classical mechanics, clustering of relativistic mechanics, etc; The methods of generalized Hamiltonian systems have been widely applied in nonlinear control theory, system control engineering, simulation algorithms and pattern recognition, especially in the fields of adaptive control, robust control, complex control, etc; In basic physical problems such as topological field theory, quantum field theory, and quantum gravity, the theory of generalized Hamiltonian systems has become an important tool for exploring field quantization due to its natural mathematical expression and physical quantitative description advantages. In summary, the research approach of the generalized Hamiltonian system reflects the overall development process of modern scientific theory, emphasizing the direction of combining mathematical tools and experiments, which will also generate value in a wider range of research fields

## 5. Generalized structural Poisson bracket

Consider smooth manifold $M$, structural function or geometric potential function $s$ characterizes its characteristics. Without loss of generality, we give the definition of generalized structural Poisson bracket in a abstract covariant form,

**Definition 5.1.** The generalized structural Poisson bracket on $\mathbb{R}^r$ of two functions $f, g \in C^\infty(M, \mathbb{R})$ is defined as

$$\{f, g\} = \{f, g\}_{GPB} + G(s, f, g)$$

where $G(s, f, g)$ is called the geobracket [2].

[2] geometric bracket is simplified as geobracket

Note that the generalized Poisson bracket is generalized by adding the geometric bracket which is directly linked to the structure of the manifold. The generalized structural Poisson bracket is formally expressed as

$$\{\cdot,\cdot\} = \{\cdot,\cdot\}_{GPB} + G(s,\cdot,\cdot)$$

This is a covariant form emerged on the manifold, the geometric bracket $G(s,\cdot,\cdot)$ is a necessary corrected term for a complete theory.

5.1. **Geometric bracket.** In this abstract representation of the generalized structural Poisson bracket, we see that geobracket equals

$$G(s,f,g) = \{f,g\} - \{f,g\}_{GPB},\ f,g \in C^{\infty}(M,\mathbb{R})$$

which always satisfies the covariant condition $G(s,f,g) \neq 0$. Actually, the geobracket is a mapping given by

$$G(s,\cdot,\cdot) : C^{\infty}(M,\mathbb{R}) \times C^{\infty}(M,\mathbb{R}) \to C^{\infty}(M,\mathbb{R})$$

$$f,g \mapsto G(s,f,g) \in C^{\infty}(M,\mathbb{R})$$

where $f,g \in C^{\infty}(M,\mathbb{R})$. As defined above, we know that there always exists the geobracket $G(s,f,g) \neq 0$ such that the generalized structural Poisson bracket holds on manifolds. Figuratively, structure function $s$ is like a stage or the background for the movement of the objects.

The non-zero of the geobracket clearly embodies the unique structure of the manifolds $M$, it's for an invariant covariance, definitely. Obviously, for this abstractly covariant form in the definition 5.1, we do not know the accurate expression of the non-zero of the geobracket, all we know is that it's a expression in terms of the structure function $s$ and other two quantities form. Note that we always ask the geobracket to satisfy the antisymmetry, that is,

$$G(s,f,g) = -G(s,g,f)$$

and one can easily verify that $\{\cdot,\cdot\}$ is skew-symmetric and bilinear. This preserves the classical properties. Meanwhile, it has some other properties such as linearity

$$G(s,\lambda f + \mu g,h) = \lambda G(s,f,h) + \mu G(s,g,h)$$

In the following discussions, we mainly discuss an accurate expression of the geobracket, maybe it's one of situation in the general covariant form. By logically derivations, we need to assure the realistic mathematical equation for the geobracket in specific form.

5.2. **Generalized gradient operator.** Firstly, we generalize the gradient operator $\nabla$ to the general generalized gradient form on manifolds.

**Definition 5.2.** The structure function $s$ related to manifolds naturally induces the generalized gradient operator to be $D = \nabla + \nabla s$ such that the generalized gradient operator of the function $f$ is

$$Df = \nabla f + f\nabla s$$

The component of generalized gradient operators is expressed as $D_i = \partial_i + \partial_i s$, called generalized derivative operators, where $\partial_i = \frac{\partial}{\partial x_i}$.

**Note that if $s = 0$ or $s =$ Constant holds, then the generalized gradient operator is degraded to the gradient operator.**

Actually, the gradient transformation is

$$\nabla \to D = \nabla + \nabla s(x)$$

equipped with the gradient of structure function $s$, note that the second part $\nabla s(x)$ of generalized gradient operator $D$ represents the character of the manifolds. For a given function $f$, we have

$$\nabla f(x) \to Df(x) = \nabla f(x) + f(x) \nabla s(x)$$

With the support of the new generalized gradient operator $D$, it directly leads to the definition of new generalized Poisson bracket as follows.

5.3. **Generalized structural Poisson bracket.** Note that $Df(x)$ is generalized gradient operator of function $f$.

**Definition 5.3** (Generalized structural Poisson bracket)**.** The generalized structural Poisson bracket is defined as

$$\{f,g\} \equiv D^T fJDg = J_{ij} D_i f D_j g$$

for $f, g \in C^\infty(M,\mathbb{R})$ on $M$, where $D = \nabla + \nabla s$ is a generalized gradient operator.

By giving the definition of generalized structural Poisson bracket, we can reasonably make sure the specific expression of geobracket $G(s,f,g)$ previously defined in 5.1.

**Theorem 5.4.** *The generalized structural Poisson bracket on* $\mathbb{R}^r$ *is written as*

$$\{f,g\} = \{f,g\}_{GPB} + G(s,f,g)$$

*for* $f, g \in C^\infty(M,\mathbb{R})$*, and geometric bracket is*

$$G(s,f,g) = f\{s,g\}_{GPB} - g\{s,f\}_{GPB}$$

*satisfying* $G(s,f,g) = -G(s,g,f)$*.*

*Proof.* Based on the definition 5.3 of generalized structural Poisson bracket and definition 5.2 of generalized gradient operator, for given $f, g \in C^\infty(M,\mathbb{R})$, it has

$$\begin{aligned}\{f,g\} &= D^T fJDg \\ &= \nabla^T fJ\nabla g + fA^T J\nabla g + g\nabla^T fJA + gfA^T JA \\ &= \{f,g\}_{GPB} + f\{s,g\}_{GPB} - g\{s,f\}_{GPB}\end{aligned}$$

where

$$\{f,g\}_{GPB} = \nabla^T fJ\nabla g,\ \{s,g\}_{GPB} = A^T J\nabla g,\ A^T J\nabla f = \{s,f\}_{GPB}$$

and $A^T JA = 0$ is given by antisymmetry. □

As explained previously, we discover a reasonable and specific expression for the abstract form of the geobracket $G(s,f,g)$ in definition 5.1, now it's sure to give this concrete description of the geobracket

$$G(s,f,g) = f\{s,g\}_{GPB} - g\{s,f\}_{GPB}$$

described by the theorem 5.4.

Hence, generalized structural Poisson bracket is written as

$$\{f,g\} = \{f,g\}_{GPB} + G(s,f,g) = \{f,g\}_{GPB} + f\{s,g\}_{GPB} - g\{s,f\}_{GPB}$$

where generalized Poisson bracket and geometric bracket both respectively satisfy Jacobi identity and rigidity theorem

$$\{f,\{g,h\}_{GPB}\}_{GPB} + \{g,\{h,f\}_{GPB}\}_{GPB} + \{h,\{f,g\}_{GPB}\}_{GPB} = 0$$
$$G(s,f,G(s,g,h)) + G(s,g,G(s,h,f)) + G(s,h,G(s,f,g)) = 0.$$

for all functions $f, g, h, s$. And, there are two other identities between geometric bracket within the generalized Poisson bracket given by

$$G(s,f,g)\{s,h\}_{GPB} + G(s,g,h)\{s,f\}_{GPB} + G(s,h,f)\{s,g\}_{GPB} = 0.$$

This inference expresses the interaction identity between geometric bracket and generalized Poisson bracket linked by geometric potential functions, which can be seen from the specific form of the expression.

$$f\{s, G(s,g,h)\}_{GPB} + g\{s, G(s,h,f)\}_{GPB} + h\{s, G(s,f,g)\}_{GPB} = 0.$$

This inference expresses the interaction identity of geometric bracket within the generalized Poisson bracket, as seen from the specific form of the expression.

Obviously, the generalized structural Poisson bracket also satisfies the antisymmetry, i.e. Namely $\{f,g\} = -\{g,f\}$. The geometric bracket $G(s,f,g)$ satisfies the antisymmetry, that is to say, $G(s,f,g) = -G(s,g,f)$. Geometric bracket also satisfies the following properties:

$$\begin{aligned} G(s,f,f) &= G(s;s,s) = 0 \\ G(s,s,g) &= s\{s,g\}_{GPB} \\ G(s,f,s) &= s\{f,s\}_{GPB} \end{aligned}$$

**Note that it is obvious that there are two cases where generalized structural Poisson bracket (GSPB) $\{f,g\} = \{f,g\}_{GPB} + G(s,f,g)$ degenerate into generalized Poisson bracket (GPB):**

- **i: If $s = 0$ or $s =$ Constant, it undergoes ordinary degradation.**
- **ii: If $G(s,f,g) = 0$, namely, geometric potential function $s = \ln(g/f)$,$fg > 0$ or $s = \ln(-g/f)$, $fg < 0$, it undergoes non-trivial degradation**

We study the generalized structure Poisson bracket with general form

$$\{f,g\} = \{f,g\}_{GPB} + f\psi + g\phi,$$

The general results are obtained under the generalized Poisson bracket framework, where

$$\psi = \psi(s,g) = \{s,g\}_{GPB}$$

and

$$\phi = \phi(s,f) = -\psi(s,f),$$

Where $s$ is a geometric potential function or structural function that is independent of $f, g$.

**Theorem 5.5.** *For $f,g,h \in C^\infty(M,\mathbb{R})$, $\lambda, \mu \in \mathbb{R}$, the generalized structure Poisson bracket have the following properties*

1. *Antisymmetry:* $\{f,g\} = -\{g,f\}$.
2. *Bilinearity:* $\{\lambda f + \mu g, h\} = \lambda\{f,h\} + \mu\{g,h\}$.
3. *Generalized Jacobi Identity:* $\{f,\{g,h\}\} + \{g,\{h,f\}\} + \{h,\{f,g\}\} = 0$.
4. *Generalized Leibnitz Identity:* $\{f,gh\} = h\{f,g\} + g\{f,h\} + gh\{s,f\}_{GPB}$.
5. *Non-degeneracy: If for $f$, $\{f,g\} = 0$ holds, then $\{f,g\}_{GPB} = G(s,g,f)$.*

*Proof.* Antisymmetry, bilinear, generalized Leibnitz identity, non-degeneracy can be easily obtained. Generalized Jacobi identity can also be obtained from the antisymmetry of generalized structure Poisson bracket. □

**Definition 5.6.** The generalized structural Poisson bracket on the smooth manifold $M$ is an operation on the smooth function space $C^\infty(M)$, for $f,g \in C^\infty(M)$, it exists $\{f,g\} \in C^\infty(M)$, the operation satisfies the condition $1-4$, $(M, S, \{,\})$ is correspondingly called generalized Poisson manifolds.

**Proposition 5.7.** *The skew-symmetric matrix $W = (W_{ij}(x))$ of order $m \times m$ is called structural matrix of the GSPB on $M$, where $W_{ij}(x) = \{x_i, x_j\}$.*

*Proof.* Based on generalized structure Poisson bracket,

$$\begin{aligned} W_{kl} = \{x_k, x_l\} &= J_{ij} D_i x_k D_j x_l \\ &= J_{ij}\left(\delta_{ik} + A_i x_k\right)\left(\delta_{jl} + A_j x_l\right) \\ &= J_{ij}\delta_{ik}\delta_{jl} + x_k J_{ij} A_i \delta_{jl} + x_l J_{ij} A_j \delta_{ik} \\ &= J_{kl} + \mathfrak{J}_{kl} \end{aligned}$$

where $A_i = \partial_i s$, $\mathfrak{J}_{kl} = x_k b_l - x_l b_k = -\mathfrak{J}_{lk}$, $b_k = J_{jk} A_j$, and $J_{kl} = J_{ij}\delta_{ik}\delta_{jl}$, $\delta_{ik}$ is Kronecker's sign. Clearly, $W_{kl} = -W_{lk}$. □

The development path of Poisson bracket can be represented by the following figure:

$$\{\cdot,\cdot\}_{CPB} \to \{\cdot,\cdot\}_{GPB} \to \{\cdot,\cdot\}_{GSPB} = \{\cdot,\cdot\}_{GPB} + G\left(s,\cdot,\cdot\right).$$

## 6. Generalized Covariant Hamilton System(GCHS)

In this section, we consider the Hamiltonian function $H$, and the generalized covariant Hamilton system is defined by generalized structure Poisson bracket.

**Theorem 6.1.** *The thorough generalized Hamiltonian systems, S-dynamics and generalized covariant Hamiltonian systems can be respectively written as*

**TGHS::** $\frac{df}{dt} = \dot{f} = \{f, H\}_{GPB} - H\{s, f\}_{GPB}$.
**S-dynamics::** $\frac{ds}{dt} = w = \{s, H\}_{GPB} = \{1, H\}$.
**GCHS::** $\frac{\mathfrak{D}f}{dt} = \{f, H\} = \{f, H\}_{GPB} + G\left(s, f, H\right) = \frac{df}{dt} + wf$.

*where* $\frac{\mathfrak{D}}{dt} = \frac{d}{dt} + w$ *is covariant time operator.*

*Proof.* According to the generalized structure Poisson bracket $\{f, g\}$, let's define generalized covariant Hamiltonian systems on $M$. Let $g = H$, then it comes to generalized covariant Hamiltonian systems.

$$\begin{aligned} \{f, H\} = D^T f J D H &= \nabla^T f J D H + f \nabla^T s J D H \\ &= J_{ij} D_i f D_j H \\ &= J_{ij} \partial_i f D_j H + f J_{ij} \partial_i s D_j H \end{aligned}$$

Let's define thorough generalized Hamiltonian systems and S-dynamics respectively

$$\begin{aligned} \frac{df}{dt} &= \nabla^T f J D H = J_{ij} \partial_i f D_j H \\ \frac{ds}{dt} &= w = \nabla^T s J D H = J_{ij} \partial_i s D_j H \end{aligned}$$

By generalized Poisson bracket $\{f, g\}_{GPB} = J_{ij} \partial_i f \partial_j g$, it easily obtains

$$J_{ij} \partial_i f D_j H = \{f, H\}_{GPB} - H\{s, f\}_{GPB}$$

$$J_{ij} \partial_i s D_j H = \{s, H\}_{GPB}$$

So we get the generalized covariant Hamilton system as desired. $\frac{\mathfrak{D}f}{dt} = \{f, H\} = \frac{df}{dt} + wf$. □

**Corollary 6.2.** *The covariant equilibrium equation is* $\frac{\mathfrak{D}f}{dt} = \{f, H\} = 0$ *such that*

$$\{f, H\}_{GPB} + G\left(s, f, H\right) = 0$$

*holds on* $M$*, it can be formally solved as* $f = f_0 e^{-wt}$*.*

In fact, the covariant equilibrium equation can be written as

$$\frac{df}{dt} = -wf \tag{6.1}$$

This equation is a typical dynamic equation. If $w \neq 0$ is a constant, then the equation (6.1) is a given homogeneous linear differential equation with constant coefficients. At this point, the properties of the solution of function $f$ are all focused on the investigation and study of S-dynamics.

In particular, let $f = s$, then the covariant equation for the structural function $s$ is

$$\frac{\mathcal{D}s}{dt} = \{s, H\} = \{s, H\}_{GPB} + G(s, s, H)$$

More specifically, it can be obtained from the generalized covariant Hamilton system.

$$\begin{aligned}\frac{\mathcal{D}s}{dt} &= \{s, H\} = \{s, H\}_{GPB} + s\{s, H\}_{GPB} \\ &= (1+s)\{s, H\}_{GPB} \\ &= (1+s)\,w\end{aligned}$$

Obviously, the covariant equilibrium equation for the structural function $s$ is $\frac{\mathcal{D}s}{dt} = 0$, namely, $(1+s)\,w = 0$, its solution is $s = -1$ or $w = 0$.

**Theorem 6.3.** *The covariant equation of the GCHS with respect to coordinates is*

$$\frac{\mathcal{D}x_k}{dt} = \{x_k, H\} = \dot{x}_k + x_k w \tag{6.2}$$

*Proof.* Based on generalized structure Poisson bracket, it easily deduces

$$\begin{aligned}\frac{\mathcal{D}x_k}{dt} &= \{x_k, H\} = J_{ij} D_i x_k D_j H \\ &= J_{ij}\partial_i x_k D_j H + x_k J_{ij}\partial_i s D_j H \\ &= J_{kj} D_j H + x_k J_{ij}\partial_i s D_j H \\ &= (J_{kj} + x_k b_j)\, D_j H\end{aligned}$$

where Kronecker sign is $\delta_{ik} = \partial_i x_k = \begin{cases}1, i = k \\ 0, i \neq k\end{cases}$. Owing to $\frac{dx_k}{dt} = J_{kj} D_j H$, thus $\frac{\mathcal{D}x_k}{dt} = \frac{dx_k}{dt} + x_k w$.

□

**Theorem 6.4.** *The TGHS is expressed in bracekt as*

$$\frac{dx_k}{dt} = \dot{x}_k = \{x_k, H\}_{GPB} + H\{x_k, s\}_{GPB}$$

Obviously, comparing the theorem 6.4 with the generalized Hamilton system (4.1), we can find that there is more one item $H\{x_k, s\}_{GPB}$ related to the structure function. This definitely is a missing part of the generalized Hamilton system. As we can see clearly, the generalized Hamilton system is just a part of thorough generalized Hamiltonian systems, it also reiterates that generalized Hamilton system is a flawed incomplete theory as introduction illustrated previously.

The global form of the generalized covariant Hamiltonian system with respect to local coordinates $x$ is written as

$$\begin{cases}\mathcal{D}x/dt = J(x)\,\nabla H(x) + H(x)\,J(x)\,\nabla s(x) + \rho(x)\,, \\ \dot{x} = J\nabla H + HJ\nabla s, \\ w(x) = \dot{s} = A^T J\nabla H,\end{cases}$$

where $\rho(x) = xw(x)$, obviously, the generalized covariant Hamiltonian system consists of three parts, namely the generalized Hamiltonian system $J(x)\nabla H(x)$and the two parts of the coupling interaction with manifold space $HJ\nabla s(x)$ and $\rho(x)$.

So the differential equation for the generalized covariant Hamiltonian system with respect to local coordinates is

$$\begin{aligned}\mathcal{D}x_j/dt = \dot{x}_i D_i x_j = \dot{x}_j + x_j w = J_{ji} D_i H + x_j w \\ = \{x_j, H\}_{GPB} - H\{s, x_j\}_{GPB} + x_j\{s, H\}_{GPB} \\ = J_{ji}\partial_i H - Hb_j + wx_j\end{aligned}$$

where

$$G(s, x_j, H) = -Hb_j + wx_j$$

is the two terms that interact or couple with the manifoldS-dynamics $w = \{s, H\}_{GPB}$.

Apparently,

**GSPB=generalized Poisson bracket+geometric bracket.**

As a more general algorithm system, Generalized Structure Poisson Brackets (GSPB) naturally defines generalized covariant Hamiltonian systems and provides many meaningful results. For Hamiltonian systems and some nonlinear systems, chaotic systems provide more mature explanations. In general, the geometric potential function is taken as a logarithmic function, such as $s(x) = \ln\zeta(x)$, where some functions with $\zeta(x) > 0$ are not constants. Therefore, it should be noted that if the geometric potential function is constant (including zero), the generalized structural Poisson bracket (GSPB) degenerates into generalized Poisson bracket (GPB), and the geometric bracket disappear, i.e. $G(s, f, g) \neq 0$. This is a special case. Therefore, in general, if the geometric potential function is not constant, i.e.$G(s, f, g) = 0$. Hence, in summary, three different algorithm systems of Poisson bracket provide corresponding Hamiltonian interpretation systems.

**Classical Poisson bracket → Classical Hamiltonian mechanics.**
**Generalized Poisson bracket (GPB) → Generalized Hamiltonian systems**
**Generalized structural Poisson bracket (GSPB)→ GCHS.**

Based on this definition of Hamilton vector field, we can define a vector field induced by the geometric potential function $s$, as follows:

**Definition 6.5.** Let $P$ be a Poisson manifold, and the geometric potential function $s : P \to \mathbb{R}$ be a smooth function, the potential vector field $X_s$ accompanying $s$ is the only smooth vector field on $P$, which satisfies

$$X_s(f) = \{f, s\}_{GPB} = -\hat{S}f$$

for all smooth functions $f : P \to \mathbb{R}$, where the structural operator $\hat{S} = b_j\partial_j = -X_s$, then $X_s$ is the geometric potential vector field on the Poisson manifold $P$.

Obviously, $X_s = c_j\partial_j = -\hat{S}$, where

$$b_j = J_{ij}A_i = -J_{ji}A_i = -c_j,$$

and $c_j = J_{ji}A_i$. From the structural operator, it can be obtained that $G(s, s, f) = -G(s, f, s) = s\hat{S}f$.

Therefore, any continuous dynamic system can find a suitable Hamiltonian energy function, and further process the dynamic system based on this energy function. We can split it into several parts such as the energy conservation part, the energy generation part, the energy consumption part, and the part that exchanges with energy. The generalized covariant

Hamiltonian system can also be represented in the form of a vector field as

$$\mathcal{D}f/dt = \{f, H\} = X_H(f) + HX_s(f) - fX_s(H)$$
$$= X_H(f) + f\hat{S}H - H\hat{S}f,$$

Where geometric bracket

$$G(s, f, H) = HX_s(f) - fX_s(H) = f\hat{S}H - H\hat{S}f,$$

is completely dependent on the potential vector field $X_s$ induced by $s$, Obviously, the generalized covariant Hamiltonian system is divided into three major force fields, namely the conservative force field and two external force fields interacting with the manifold.

Therefore, the generalized structural Poisson bracket can also be expressed as

$$\{f, g\} = X_g f + gX_s f - fX_s g,$$

where $X_s = c_j\partial_j$ is the geometric potential vector field on $M$, which is obviously easy to obtain $X_s x_j = c_j$, and $X_g = -J_{ij}\partial_i g\partial_j$, where $X_g = -J_{ij}\partial_i g\partial_j$ is the vector field induced by the smooth function $g$ on the Poisson manifold $P$. The geometric bracket can also be represented as

$$G(s, f, g) = f\hat{S}g - g\hat{S}f,$$

Where $\hat{S}$ is the structural operator on $M$ and $c_j = -b_j$, where

$$X_s f = \{f, s\}_{GPB} = -\{s, f\}_{GPB} = -\hat{S}f = -X_f s,$$

And $\hat{S}f = X_f s = \{s, f\}_{GPB}$ is the derivative with respect to the geometric potential function $s$ in the direction of a vector field $X_f$. Like S-dynamics, $w = \hat{S}H = X_H s$ is the derivative of the direction along the Hamiltonian vector field $X_H$ of the geometric potential $s$.

The difference between generalized structural Poisson bracket and classical generalized Poisson bracket is mainly reflected in geometric bracket. As a parallel existence with classical generalized Poisson brackets, they jointly act on a more general and extensive Hamiltonian system, including chaotic systems, nonlinear systems, and other systems. Geometric brackets serve as a beneficial supplement to classical generalized Poisson brackets, It has an independent theoretical structure system similar to the classical generalized Poisson bracket, which can be seen from the research results of geometric bracket.

The development prospects of generalized structural Poisson bracket is very broad. Firstly, the generalized Poisson bracket can delve deeper into gravitational problems, providing an effective mathematical tool for understanding these problems. On the other hand, in mathematical research, the generalized structural Poisson bracket also have great potential. For example, it can help delve deeper into mathematical problems such as Lie algebra and Hamilton Jacobi equations. The generalized structure Poisson bracket can also be extended to applications in semi-classical and geometric quantum mechanics, which will further expand the application and understanding of quantum mechanics under the generalized covariant limit. In short, in the future, the generalized structure Poisson bracket will continue to play an important role in providing an important mathematical tool for theoretical research and applications in fields such as general relativity and physical mathematics.

## 7. Conclusions

The generalized Hamilton system on $\mathbb{R}^r$ is a natural extension of the classical Hamilton system on $\mathbb{R}^{2n}$. The corresponding classical Poisson brackets is extended to the generalized Poisson brackets, which promotes the unprecedented development of Hamilton mechanics. As discussed in the introduction, it's extended to any dimension. but generalized Hamilton system can not deal with the differential equation of fixed-point rotation of rigid body under external moment, some practical problems can not be expressed by generalized Hamilton system, Hamiltonian realization problem, the basic theoretical framework of generalized Hamilton system is not perfect, and many other problems. The generalized covariant Hamilton system is a complete theoretical system on $\mathbb{R}^r$. As a covariant extension of the generalized Hamilton system, it contains two independent systems: the generalized thorough Hamiltonian system and the S-dynamic system. The latter, as an independent system, is only related to the structure function, which is induced by the structure function and mainly describes the rotation problem. Thus, there is no needing to add an additional term artificially, a complete generalized covariant Hamilton system automatically includes a patch term that describes the fixed-point rotation problem. Generalized covariant Hamilton system defined by generalized structure Poisson bracket shows that the dynamic evolution of function $f$ is related to the specific structure function of manifold and is constrained by manifold structure. Generalized covariant Hamilton system is constrained by manifold structure because of the perfection of structure function. It can always be expressed as a nonlinear system. From the equilibrium solution of the generalized covariant Hamilton system, it can be seen that it is a nonlinear solution, and a general dynamic solution. Obviously, the dynamic steady state property of the system is completely evolved by the evolution of S-dynamic system over time. The introduction of geometric brackets is necessary for generalized Poisson bracket after generalization. It has successfully linked the dynamic evolution of manifold structures and functions, greatly expanded the scope of research.

(Gen WANG)
Department of Mathematics,
ZheJiang Normal University, ZheJiang, JinHua 321004,
P.R.China
*Email address*: `wanggen@zjnu.edu.cn`